\documentclass{article}
\usepackage{amssymb}
\usepackage{amsmath}

\newtheorem{theorem}{Theorem}[section]

\newtheorem{corollary}[theorem]{Corollary}

\newtheorem{lemma}[theorem]{Lemma}

\newtheorem{proposition}[theorem]{Proposition}

\newenvironment{proof}[1][Proof]{\noindent\textbf{#1.} }{\ \rule{0.5em}{0.5em}}

\begin{document}

\title{\bf \Large Convergence of fundamental solutions of linear parabolic
equations under Cheeger-Gromov convergence}

\author{Peng Lu}

\date{May 03, 2010}
%%%relative definite version: May 03, 2010
\maketitle

\section{Introduction}
Heat kernel is an important tool and has many applications in differential 
geometry (see, for example, \cite{LY}) and physics.
Heat kernel and related ideas play an important role in Perelman's work 
on Ricci flow, one particular example is the so-called Perelman's Li-Yau type 
differential Harnack inequality(\cite[Corollary 9.3]{Pe02I}, see also \cite{Ni}
and \cite{CTY} for more details of the proof).
  
In the proof of the pseudo locality
theorem (\S 10.1, in particular, p.25 in \cite{Pe02I})
Perelman used the following statement about the convergence of heat kernels.
Let $(M^n_k,g_k(\tau),x_{0k})$, $\tau \in [0,T]$, be a pointed sequence
of $n$-dimensional complete solutions of the backward Ricci flow.
Suppose that each solution has bounded curvature and suppose that the sequence
converges to $(M^n_\infty, g_\infty(\tau), x_{0\infty})$ in the pointed
Cheeger-Gromov sense, then the heat kernels $H_k$ for the
adjoint heat equations associated with $(M_k,g_k(\tau),x_{0k})$
sub-converge to a positive limit $H_\infty$ defined on $M_\infty \times(0,T]$.
Motivated by this, the proof of the statement is studied by several
groups of authors, Chau, Tam and Yu(\cite{CTY}), Hsu(\cite{Hs}),
Kleiner and Lott(\cite{KL}), S. Zhang(\cite{ZhS}), and
Chow(\S 22.2 in \cite{CCIII}).

In this note by combining the techniques from the literature mentioned above
we give a proof of the convergence of fundamental solutions for general
linear parabolic equations(see Theorem \ref{thm de main conv} below).
More precisely our proof of the existence of the limit is similar
to the one given by Chow in \cite{CCIII}.
Our proof of the $\delta$-function property of the limit is very close to
the proofs given in \cite{CTY} and \cite{Hs}.
Another feature worth mention is that we do not assume that the sequence has
uniformly bounded curvature (compare to \cite[Theorem 0.1]{ZhS}).
This is made possible by avoiding the use of decay estimates of heat kernels
in our proof.

The above mentioned proof of Theorem \ref{thm de main conv} will be given in \S 2.
In \S 3 we will discuss when the assumption, about the upper bound of $L^1$-norms
made in Theorem \ref{thm de main conv}, could be satisfied.
In \S 4 we consider the uniqueness of the limit given by Theorem
\ref{thm de main conv}.
In \S 5 we construct a cut-off function following Perelman (\cite[\S 8.3]{Pe02I})
and use it to prove a local
integral estimate of fundamental solutions(see Corollary
\ref{cor local int fund sol lower bdd}),
essentially this is the only novel part in this note.

%%%%%%%%%%%%%%%%%%%%%%%%%%%%%%%%%%%%%%%%%%%%%%%%%%%%%%%%%%%%%%%%%%%%%%%%%%%%
%%%%%%%%%%%%%%%%%%%%%%%%%%%%%%%%%%%%%%%%%%%%%%%%%%%%%%%%%%%%%%%%%%%%%%%%%%%%
\section{The convergence of fundamental solutions}

\textbf{2.1 The fundamental solutions}.
Let $M^n$ be a smooth connected manifold and
let $\Omega \subset M$ be a smooth connected
domain with (possibly empty) smooth boundary $\partial\Omega$.
Fix a $T>0$, let $g\left(\tau\right)$, $\tau\in\left[  0,T\right]$, be
a smooth family of smooth metrics on $\Omega$.
Define the symmetric $2$-tensor $\mathcal{R}_{ij}$
and the scalar function $\mathcal{R}$ by
\begin{align}
& \mathcal{R}_{ij}\doteqdot\frac{1}{2}\, \frac{\partial g_{ij}}{\partial\tau},
\label{evol eq gen} \\
& \mathcal{R} \left(  x,\tau\right)  \doteqdot g^{ij}\left(  x,\tau\right)
\mathcal{R}_{ij}\left(  x,\tau\right)  . \label{eq tr ricci gen}
\end{align}

Let $Q: \Omega \times\left[  0,T\right] \rightarrow\mathbb{R} $
be a smooth function.
Let $X(\tau)$, $\tau \in [0,T]$, be a smooth family of smooth
vector fields on $\Omega$.
We will consider the following linear parabolic equation on $\Omega
\times [0,T]$ with respect to the evolving metrics,
\begin{equation}
\square^{\ast}u \doteqdot \frac{\partial u}{\partial\tau}
- \Delta_{g\left(  \tau\right)  }u + \nabla_{X}
 u + Q u =0 \label{eq heat gen}
\end{equation}
where $\nabla_{X}u \doteqdot X(\tau)u$.

Let $\operatorname{Int}(\Omega)$ denote the set of interior points of $\Omega$.
Fix a base point $x_{0}\in \operatorname{Int}(\Omega)$.
Let $H\left( x,\tau \right)$ be a function belong to
\[
C^0( \Omega \times (0,T] , \mathbb{R}) \cap C^\infty(
 \operatorname{Int}(\Omega) \times(0,T] , \mathbb{R}) .
 \]
$H$ is called a \textbf{fundamental solution} of
(\ref{eq heat gen}) centered at $x_0$ if it satisfies the following.
$H$ is positive on $\operatorname{Int}(\Omega) \times(0,T]$ and
\begin{eqnarray}
\square^{\ast} H &=&0  \label{eq fund sol} \\
\lim_{\tau \rightarrow 0^+} H \left(\cdot, \tau \right)  &=&\delta _{x_{0}}.  \notag
\end{eqnarray}
$H$ is called the \textbf{heat kernel} of
(\ref{eq heat gen}) centered at $x_0$ if $H$ is the minimal fundamental solution
of (\ref{eq heat gen}) centered at $x_0$.

\noindent \textbf{Remark}. Fundamental solutions as defined above are not necessarily
unique. In particular when $\Omega$ is a compact manifold with nonempty boundary.
Let $h$ be a solution of $\square^{\ast} h =0$ and which satisfies boundary conditions
$\left . h \right |_{\Omega \times \{0 \}} =0$ and
$\left . h \right |_{\partial \Omega \times [0,T]} \geq 0$.
If $H$ is a fundamental solution, then $H+h$ is also a fundamental solution.

\noindent \textbf{Remark}. In this note we do not address the issue of the existence
of $H_k$. It is well-known that convergence theorems, such as Theorem
\ref{thm de main conv}, can be used to prove the existence of heat kernels on
complete noncompact manifold assuming the existence of Dirichlet heat
kernels on compact manifolds with boundary.

%%%%%%%%%%%%%%%%%%%%%%%%%%%%%%%%%%%%%%%%%%%%%%%%%%%%%%%%%%%%%%%%%%%%%%%%%%%%
%%%%%%%%%%%%%%%%%%%%%%%%%%%%%%%%%%%%%%%%%%%%%%%%%%%%%%%%%%%%%%%%%%%%%%%%%%%%
\vskip 0.2cm
\noindent \textbf{2.2 The main theorem}.
Let $\left \{ M_k^{n} \right \}_{k=1}^\infty$ be a sequence of smooth manifolds,
let $\Omega_{k} \subset M_k $ be a smooth connected domain with (possibly empty)
smooth boundary $\partial \Omega_k$ for each $k$,
and let $g_{k}\left(  \tau \right)$, $\tau \in\left[0, T \right]$,
be a smooth family of smooth metrics on $\Omega_k$ for each $k$.
Let $x_{0k} \in \operatorname{Int}(\Omega_k)$.
We assume that the sequence $\left \{ \left (\Omega_k,g_k(\tau), (x_{0k},0) \right)
\right \}$ converges, in the $C^{\infty}$ pointed Cheeger-Gromov sense,
to a smooth limit $\left( M_{\infty}^n,g_{\infty}\left(  \tau\right),
\left( x_{0\infty}, 0\right)  \right)$, $\tau\in\left[  0, T \right]$.
Here $ \left( M_{\infty}, g_{\infty}(\tau) \right)$ may have boundary
and may not be complete.
By definition the convergence means that there exist an exhaustion
$\left\{  U_k \right \}_{k=1}^ \infty$ of $M_{\infty}$ by
pre-compact open sets and a sequence of diffeomorphisms
$ \Phi_k: U_k \rightarrow  V_k
\doteqdot \Phi_{k}\left(  U_k \right) \subset \Omega_k$
such that \\
(\textbf{CG1}) $x_{0\infty} \in U_k$ and $\Phi_k \left(x_{0\infty}\right)=x_{0k}$
for each $k$, and \\
(\textbf{CG2}) $\left( U_k,\Phi_k^{\ast}\left[  \left. g_k \left(
\tau\right) \right\vert_{ V_k }\right]  \right)$
converges uniformly in $C^{\infty}$-topology to $g_{\infty} \left(\tau\right)$
on any compact subset of $ M_{\infty}\times\left[  0, T\right]$. \\
From the definition it is easy to see that $x_{0\infty} \in \operatorname{Int}(M_\infty)$.

Let $Q_k: \Omega_k \times\left[  0,T\right] \rightarrow\mathbb{R} $
be a smooth function for each $k$.
Let $X_k(\tau)$, $\tau \in [0,T]$, be a smooth family of smooth vector fields
on $\Omega_k$ for each $k$.
Suppose that under the convergence of $\left \{ \left (\Omega_k,g_k(\tau),
(x_{0k},0) \right) \right \}$ the sequence $\{Q_k \}$ converges to a smooth limit
$Q_\infty: M_\infty \times \left[0,T\right] \rightarrow \mathbb{R}$,
i.e., the sequence of functions $\{Q_k (\Phi_k (\cdot ), \tau)  \}$
converges uniformly in $C^{\infty}$-topology to $Q_{\infty} \left(\cdot, \tau\right)$
on any compact subset of $ M_{\infty}\times \left[0, T\right]$.
We also assume that the sequence of vector fields $\{X_k(\tau) \}$
converges to a smooth vector field limit $X_\infty( \tau) $
on $ M_\infty \times \left[0,T\right]$,
i.e., the sequence of push-forward vector fields $\{ \left (\Phi_k^{-1}
\right )_* X_k( \tau)   \}$
converges uniformly in $C^{\infty}$-topology to $X_{\infty} \left( \tau\right)$
on any compact subset of $ M_{\infty}\times \left[0, T\right]$.

We will consider the following equation on $\Omega_k \times [0,T]$.
\begin{equation}
\square_k^{\ast} u \doteqdot \frac{\partial u}{\partial \tau }
-\Delta _{g_k \left( \tau \right) }u + \nabla_{X_k}
 u +Q_k u =0.   \label{eq seq fund sol}
\end{equation}
Let $H_{k}:\Omega_k \times(0, T]\rightarrow\left(  0,\infty\right)
$ be a fundamental solution of (\ref{eq seq fund sol}) centered at $x_{0k}$,
i.e., $\square_k^{\ast} H_k =0$ and
\begin{equation}
\lim_{\tau \rightarrow 0} H_k(\cdot,\tau) =\delta_{x_{0k}}. \label{eq Hk is delta k}
\end{equation}
We define function $f_k: \Omega_k \times(0, T]\rightarrow \mathbb{R}$ by
\begin{equation}
H_k (x,\tau) \doteqdot (4 \pi \tau)^{-n/2} e^{-f_k (x,\tau)}.
\label{eq f in seq fund sol}
\end{equation}

The following is the convergence property of fundamental solutions,
under the pointed $C^\infty$ Cheeger-Gromov convergence of the underlying manifolds,
the convergence of potential functions $Q_k$, and the convergence of vector
fields $X_k$.

\begin{theorem} \label{thm de main conv}
Under the setup and notations above,
suppose there is a constant $C_\ast$ independent of $k$ such that
\begin{equation}
\int_{\Omega_k }
H_{k}\left( x,\tau \right) d\mu _{g_{k}\left( \tau \right)}(x)
\leq C_\ast   \label{eq Hk int assump}
\end{equation}
for each $k$ and $\tau \in ( 0,T]$,
then there is a subsequence (still indexed by $k$)
such that the following holds.
There are smooth nonnegative function $H_{\infty}$ defined on
$\operatorname{Int}(M_\infty) \times ( 0,T] $ and smooth function $f_{\infty}$
defined on $\operatorname{Int}(M_\infty) \times (0,T]$ satisfying
\begin{equation}
H_\infty (x,\tau) = (4 \pi \tau)^{-n/2} e^{-f_\infty(x,\tau)}
\label{eq H f relation at infty},
\end{equation}
such that
\begin{equation}
\tilde{H}_k  (\cdot,\cdot) \doteqdot H_k (\Phi_k(\cdot), \cdot)
\longrightarrow H_{\infty} (\cdot,\cdot) \label{eq H converges}
\end{equation}
uniformly in $C^{\infty}$-topology on any compact subset of
$\operatorname{Int}( M_\infty) \times (0,T] $, and
\begin{equation}
\tilde{f}_k  (\cdot,\cdot)  \doteqdot f_k (\Phi_k (\cdot), \cdot)
\longrightarrow f_{\infty} (\cdot,\cdot) \label{eq f converges}
\end{equation}
uniformly in $C^{\infty}$-topology on any compact subset of
$\operatorname{Int}(M_\infty) \times (0,T]$.
Moreover, $H_{\infty}$ is a fundamental solution to the parabolic
equation on $\operatorname{Int}(M_\infty) \times [0,T]$
\begin{align}
& \square_{\infty}^{\ast}H_{\infty}\doteqdot\left(  \frac{\partial}{\partial
\tau}-\Delta_{g_{\infty(\tau)}} + \nabla_{X_\infty}
+Q_{\infty}\right)  H_{\infty}=0 \label{eq Hinfty sat} \\
&  \lim_{\tau \rightarrow 0^+} H_\infty \left(\cdot, \tau \right)
=\delta _{x_{0\infty}}, \label{eq delta thm main 1}
\end{align}
and
\begin{equation}
\int_{M_\infty }
H_\infty \left( x,\tau \right) d\mu _{g_\infty \left( \tau \right)}(x)
\leq C_\ast   \label{eq H infty upper bound finite}
\end{equation}
for each $\tau \in (0,T]$.
\end{theorem}

\noindent \textbf{Remark}. From the proof(see (\ref{eq what is what use 01}) and
(\ref{eq uk F ode bdd})),
the assumption (\ref{eq Hk int assump}) can be replaced by the following
condition and we still have the convergence to a fundamental solution.
For any compact domain $\hat{D} \subset \operatorname{Int}(M_\infty)$
there is a constant $C_{\hat{D}}$
independent of $k$ such that when $k$ large enough
\begin{equation}
\int_{\Phi_k(\hat{D}) }
H_{k}\left( x,\tau \right) d\mu _{g_{k}\left( \tau \right)} (x)
\leq C_{\hat{D}}   \label{eq Hk int assump mod}
\end{equation}
for any $\tau \in (0,T]$. However under assumption (\ref{eq Hk int assump mod})
the limit $H_\infty$ may not satisfy (\ref{eq H infty upper bound finite}).

\noindent \textbf{Remark}. When $M_\infty$ in Theorem \ref{thm de main conv}
is a compact manifold with nonempty boundary, the theorem does not address
the issue about whether the convergence can be extended to the boundary.

\vskip .1cm
The remaining of this section is devoted to give a proof of
Theorem \ref{thm de main conv}.

%%%%%%%%%%%%%%%%%%%%%%%%%%%%%%%%%%
\vskip 0.2cm
\noindent \textbf{2.3 The mean value inequality}.
The following parabolic mean value inequality will be used in the proof of
Theorem \ref{thm de main conv}.

Let $\Omega \subset M^n$ be a smooth connected domain in a smooth manifold.
Let $g(\tau))$, $\tau \in [0,T]$, be a smooth family of smooth
metrics on $\Omega$.  We adopt the notations
$\mathcal{R}_{ij}$, $\mathcal{R}$, $X$, and $Q$ from the beginning of \S2.1.

Let $\tilde{g}$ be a smooth metric on $\Omega$
satisfying $\operatorname{Rc}\left( \tilde{g}\right) \geq -K$ on $\Omega$
for some $K \geq 0$. We assume that
\begin{equation}
C_0^{-1} \tilde{g} \leq g(0) \leq C_0 \tilde{g} \label{eq g tilde g0 equ}
\end{equation}
in $\Omega$ for some constant $C_0$.

\begin{theorem} \label{thm mvi}
Under the setup and the assumption above.
Let $u:\Omega\times\left[  0,T\right] \rightarrow\mathbb{R}$ be a nonnegative
sub-solution to (\ref{eq heat gen}), i.e.,
\[
\frac{\partial u}{\partial\tau} - \Delta_{g\left(  \tau\right) }u
+ \nabla_Xu+ Q\,u  \leq 0 .
\]
If $\left( x_{0},\tau_{0}\right)  \in\Omega\times(0,T]$ and $r_{0}>0$ are such that
the closed ball $\overline{B_{\tilde{g}}\left(  x_{0}, 2r_0 \right) }$ is
contained in $\operatorname{Int}(\Omega)$ and is compact.
We assume $\tau_0  \geq (2r_0)^2$ so that the parabolic cylinder
\begin{equation}
\overline{B_{\tilde{g}}\left(  x_{0}, 2 r_0 \right) } \times
[\tau_{0}-\left(  2r_{0}\right)
^{2}, \tau_0 ] \subset \Omega\times\left[  0,T\right]  .
\label{No More Footnotes}
\end{equation}
Then for any $(x, \tau ) \in B_{\tilde{g}}\left(  x_{0}, r_0 \right)  \times [\tau_{0}-r_{0}^{2}, \tau_{0}]$ we have
\begin{equation}
u (x, \tau) \leq\frac{C_{1}e^{C_{2}\tau_{0}+C_{3}\sqrt{K}r_{0}}}{r_{0}^{2}%
\operatorname{Vol}_{\tilde{g}}\left(  B_{\tilde{g}}\left(  x_{0},r_{0}\right)
\right)  }\int_{ \overline{B_{\tilde{g}}\left(  x_{0}, 2r_0 \right) } \times
[\tau_{0}-\left(  2r_{0}\right)^{2}, \tau_0 ] }\!u\!\left( y, \sigma \right)  d\mu_{\tilde{g} }\!\left( y \right)  d\sigma, \label{eq mvi}
\end{equation}
where constant $C_{1}$ depends only on $n$, $T$, $C_{0}$, and $\sup_{\Omega\times
\left[  0,T\right]  }\left\vert \mathcal{R}_{ij}\right\vert _{g\left(  \tau\right)}$,
constant $C_{2}$ depends only on $\sup_{\Omega\times\left[  0,T\right]
} \max \{-Q,0 \}$, $\sup_{\Omega\times\left[  0,T\right]
}\left\vert X \right\vert_{g(\tau)}$ and $\sup_{\Omega\times\left[  0,T\right]
}\left\vert \mathcal{R}  \right\vert$,
and constant $C_{3}$ depends only on $n$.
\end{theorem}

\noindent \textbf{Remark} When $X=0$ Theorem \ref{thm mvi} is Lemma 3.1 in
\cite{CTY}, which is based on \cite[\S 5]{ZhQ}.
The following proof of the mean value inequality is a slight modification
of the proof given in \cite[\S25.1]{CCIII} which also considers the case $X=0$.

\vskip .1cm
\noindent \textbf{Proof}
Given a nonnegative sub-solution $u:\Omega\times\left[  0,T\right]  \rightarrow\mathbb{R}_{+}$ to (\ref{eq heat gen}), define
\begin{equation}
v\doteqdot e^{-A\tau}u, \label{v definition}
\end{equation}
where $A \geq 0$ to be chosen later (see (\ref{eq A choice}) below).
We have
\[
\frac{\partial v}{\partial\tau}-\Delta v
+\nabla_{X}v +\left(  Q+A\right)  v \leq 0
\]
where we have dropped ${g(\tau)}$ from our notation, e.g.,
$\Delta=\Delta_{g(\tau)}$.
Hence for any real number $p\in\lbrack1,\infty)$
\begin{equation}
\frac{\partial}{\partial\tau}\left(  v^{p}\right)  -\Delta\left(
v^{p}\right) +\nabla_X v^p +p\left(  Q+A\right)  v^{p}\leq-p\left(  p-1\right)
v^{p-2}\left\vert \nabla v\right\vert ^{2}\leq 0 \label{eq vp inequ MV}
\end{equation}
on $\Omega\times\left[  0,T\right]  $.
To prove the so-called reverse Poincare-type inequality (see
(\ref{eq reverse sobolev ineq}) below),
we will first localize this inequality and then integrate it.

Let $0\leq\tau_{1}<\tau_{2}\leq T$ and let%
\begin{equation}
\psi:\Omega\times\left[  \tau_{1},\tau_{2}\right]  \rightarrow\left[
0,1\right]  \label{eq cutoff psi}
\end{equation}
be a cutoff function with support contained in $D\times\left[  \tau_{1}
,\tau_{2}\right]  $, where $D\subset\Omega$ is a compact domain with
$C^1$-boundary. Assume that
\begin{equation}
\psi\left(  x,\tau_{1}\right)  \equiv 0  \,\,\,\,
\text{ for } x\in \Omega. \label{eq psi initial 0}
\end{equation}
Furthermore we assume $\frac{\partial\psi}{\partial\tau} \geq 0$ and
\begin{equation}
\psi\frac{\partial\psi}{\partial\tau}+\left\vert \nabla\psi\right\vert
_{g\left(  \tau\right)  }^{2}\leq L \quad \text{ on }
D\times\left[ \tau_{1},\tau_{2}\right] \label{eq psi diff ineq}
\end{equation}
for some constant $L\in\lbrack0,\infty)$.

Multiplying the inequality (\ref{eq vp inequ MV}) by
$\psi^{2}v^{p}$ and integrating by parts in space and time, we have
\begin{align}
0  &  \geq\int_{\tau_{1}}^{\tau_{2}}\int_{D}\psi^{2}\left(  \frac{1}{2}
\frac{\partial}{\partial\tau}\left(  v^{2p}\right)  -v^{p}\Delta\left(
v^{p}\right) +  v^p \nabla_X v^p  +p\left(  Q+A\right)
v^{2p}\right)  d\mu\,d\tau \notag \\
&  =\int_{\tau_{1}}^{\tau_{2}}\int_{D}\left(  -\psi\frac{\partial\psi
}{\partial\tau}v^{2p}-\frac{1}{2}\psi^{2}v^{2p}\mathcal{R}\right)  d\mu
\,d\tau+\frac{1}{2}\left(  \int_{D}\psi^{2}v^{2p}d\mu\right)  \left(  \tau
_{2}\right)  \notag \\
&  \quad\;+\int_{\tau_{1}}^{\tau_{2}}\int_{D}\left(  \left\vert \nabla\left(
\psi v^{p}\right)  \right\vert ^{2}-v^{2p}\left\vert \nabla\psi\right\vert
^{2}+p\left(  Q+A\right)  \psi^{2}v^{2p}\right)  d\mu\,d\tau \notag \\
& \quad\; +\int_{\tau_{1}}^{\tau_{2}}\int_{D}\left( \psi v^p \nabla_X (\psi v^p)
- \psi v^{2p} \nabla_X \psi  \right)  d\mu\,d\tau \label{eq vp must be taged}
\end{align}
where we have used $\frac{\partial}{\partial\tau}d\mu =\mathcal{R\,}d\mu$ and
\[
-\int_{D}\psi^{2}v^{p}\Delta\left(  v^{p}\right)  d\mu=\int_{D}\left\vert
\nabla\left(  \psi v^{p}\right)  \right\vert ^{2}d\mu-\int_{D}v^{2p}\left\vert
\nabla\psi\right\vert ^{2}d\mu
\]
to derive the equality.

Let $A_1 \doteqdot 1+ \sup_{\Omega\times\left[  0,T\right]
}\left\vert X \right\vert_{g(\tau)}$. We compute
\begin{align}
& \left \vert \int_{\tau_{1}}^{\tau_{2}}\int_{D} \psi v^p \nabla_X (\psi v^p)
 d\mu\,d\tau \right \vert \notag \\
 \leq \, &  A^2_1 \int_{\tau_{1}}^{\tau_{2}}\int_{D} \psi^2 v^{2p}
 d\mu\,d\tau + \frac{1}{4A^2_1} \int_{\tau_{1}}^{\tau_{2}}\int_{D} \vert \nabla_X
 (\psi v^p) \vert^2 d\mu\,d\tau \notag \\
 \leq \, &   A^2_1 \int_{\tau_{1}}^{\tau_{2}}\int_{D} \psi^2 v^{2p}
 d\mu\,d\tau + \frac{1}{4} \int_{\tau_{1}}^{\tau_{2}}\int_{D} \vert \nabla
 (\psi v^p) \vert^2 d\mu\,d\tau. \label{eq mod chow 1}
\end{align}
We estimate
\begin{align}
& \left \vert \int_{\tau_{1}}^{\tau_{2}}\int_{D} \psi v^{2p} \nabla_X \psi
d\mu\,d\tau \right  \vert \notag \\
\leq \, & \frac{1}{4}A_1^2 \int_{\tau_{1}}^{\tau_{2}}\int_{D} \psi^2 v^{2p}
 d\mu\,d\tau + \frac{1}{A_1^2} \int_{\tau_{1}}^{\tau_{2}}\int_{D}
  v^{2p} \vert  \nabla_X \psi \vert^2 d\mu\,d\tau \notag \\
  \leq \, &  \frac{1}{4}A_1^2 \int_{\tau_{1}}^{\tau_{2}}\int_{D} \psi^2 v^{2p}
 d\mu\,d\tau + \int_{\tau_{1}}^{\tau_{2}}\int_{D}
  v^{2p} \vert  \nabla \psi \vert^2 d\mu\,d\tau. \label{eq mod chow 2}
\end{align}
Hence by combining (\ref{eq vp must be taged}), (\ref{eq mod chow 1}),
and (\ref{eq mod chow 2}) we get
\begin{align}
0  &  \geq \frac{3}{4} \int_{\tau_{1}}^{\tau_{2}}\int_{D}\left\vert
\nabla\left(  \psi v^{p}\right)  \right\vert ^{2}d\mu\,
d\tau+\frac{1}{2}\left(  \int_{D}\psi
^{2}v^{2p}d\mu\right)  \left(  \tau_{2}\right) \notag  \\
&  \quad\;-\int_{\tau_{1}}^{\tau_{2}}\int_{D}\left(  \psi\frac{\partial\psi
}{\partial\tau}+ 2 \left\vert \nabla\psi\right\vert ^{2}\right)  v^{2p}
d\mu\,d\tau  \notag  \\
&  \quad\; + \int_{\tau_{1}}^{\tau_{2}}\int_{D} \left(p(Q+A) -\frac{1}{2}
\mathcal{R} - \frac{5}{4}A_1^2 \right) \psi^2 v^{2p}
 d\mu\,d\tau .  \label{eq 26a add}
\end{align}

Now choose $A\in\lbrack0,\infty)$, depending only on $\sup
_{\Omega\times\left[  0,T\right]  } \max \{ -Q(x,\tau),0\} $, $\sup
_{\Omega\times\left[  0,T\right]  }\mathcal{R}$, and $\sup_{\Omega\times\left[  0,T\right]}\left\vert X \right\vert_{g(\tau)}$(in particular, $A$ is
independent of $p$ and $D$), so that
\begin{equation}
Q+A-\frac{1}{2p}\mathcal{R} - \frac{5}{4p}A_1^2 \geq 0 \label{eq A choice}
\end{equation}
on $D\times\left[  0,T\right]  $ for all $p \geq 1$.
Hence from (\ref{eq psi diff ineq})
and (\ref{eq 26a add}) we have the following.

\begin{lemma} \label{lem reverse Sobolev}
For the choice of $A$ as in (\ref{eq A choice}) the function $v$ defined
in (\ref{v definition}) satisfies
\begin{equation*}
 \int_{\tau_{1}}^{\tau_{2}}\int_{D}\left\vert \nabla\left(  \psi v^{p}\right)
\right\vert _{g\left(  \tau\right)  }^{2}d\mu_{g\left(  \tau\right)  }%
\,d\tau\leq \frac{8}{3}L\int_{\tau_{1}}^{\tau_{2}}\int_{D}v^{2p}d\mu_{g\left(
\tau\right)  }\,d\tau
\end{equation*}
and
\begin{equation}
 \left(  \int_{D}\psi^{2}v^{2p}d\mu_{g\left(  \tau\right)  }\right)  \left(
\tau_{2}\right)  \leq 4 L\int_{\tau_{1}}^{\tau_{2}}\int_{D}v^{2p}d\mu_{g\left(
\tau\right)  }\,d\tau, \label{eq reverse sobolev ineq}
\end{equation}
where $D\subset\Omega$ is a compact domain with $C^1$-boundary and where $\psi$
satisfies $(\ref{eq psi diff ineq})$  and  $\operatorname*{supp}\left(  \psi\right)  \subset D\times\left[ \tau_{1},\tau_{2}\right]  $.
\end{lemma}

The remaining proof of Theorem \ref{thm mvi} is the same as the proof
given in \S25.1 of \cite{CCIII}, we omit it.
Now we turn to the proof of Theorem \ref{thm de main conv}.

%%%%%%%%%%%%%%%%%%%%%%%%%%%%%%%%%%%%%%%%%%%%%%%%%%%%%%%%%%%%%%%%%%%%%%%%%%%%
%%%%%%%%%%%%%%%%%%%%%%%%%%%%%%%%%%%%%%%%%%%%%%%%%%%%%%%%%%%%%%%%%%%%%%%%%%%%
%%%  \newpage
\vskip .2cm
\noindent \textbf{2.4 Step 1 The existence of $H_\infty$}. First we apply Theorem
\ref{thm mvi} to $H_k$ in Theorem \ref{thm de main conv} to get some local
$C^0$ estimate of $H_k$ uniform in $k$.
In this subsection we adopt the notation of \S2.2.
Fix an arbitrary interval $\left[ \tau _{1},\tau _{2}\right] \subset (0,T]$ and
fix a compact domain set $D \subset \operatorname{Int}(M_\infty)$
which contains $x_{0\infty}$.
Let $\tilde{g}=g_{\infty }\left( 0\right) $. We define for any $r >0$
\[
D \left( r \right) \doteqdot \cup _{x \in D} B_{\tilde{g}}\left( x,r\right)
\subset M_\infty.
\]
Below we fix a constant $r_1>0$ such that the closure $\overline{D \left( r_1
\right)}$ is compact in $\operatorname{Int}(M_\infty)$.

Choose $k_0$ large enough such that $D(r_1) \subset U_k$ for all $k \geq k_0$.
Since $g_k(\tau)$, $X_k$, and $Q_k$,
after adjustments by $\Phi_k$, converge to $g_\infty (\tau)$,
 $X_\infty$, and $Q_\infty$ respectively and uniformly in $C^\infty$-norm
on any compact subset of  $\operatorname{Int}(M_\infty)$,
there is a $C_0$ independent of $k$ but depending on $D \left( r_1 \right)$
such that for any $k \geq k_0$
\begin{align}
& C_0^{-1} \tilde{g}\leq \Phi_k^* g_k \left( \tau \right) \leq C_0 \tilde{g} \,\,\,\,\,
\text{ on } D\left( r_1 \right)  \text{ for all } \tau \in [0,T],
\notag \\
& \sup_{\Phi_k  \left( D(r_1) \right) \times \left[ 0,T\right] }
\left\vert \left( \mathcal{R}
_{k}\right) _{ij} \right\vert _{g_{k}\left( \tau
\right) }\leq  C_0 \left (\sup_{D \left( r_1 \right) \times \left[ 0,T\right]
}\left\vert \left( \mathcal{R}_{\infty }\right) _{ij}
\right\vert _{g_{\infty }\left( \tau \right) } +1 \right) , \notag \\
& \sup_{\Phi_k \left(D( r_1) \right) \times \left[ 0,T \right] } \left\vert \mathcal{R}_{k} \right\vert \leq C_0 \left ( \sup_{D
\left( r_1\right) \times
\left[ 0,T\right] }\left\vert \mathcal{R}_{\infty }
\right\vert +1 \right ) , \label{eq get k rid of} \\
& \sup_{\Phi_k \left(D( r_1) \right) \times \left[ 0,T \right] }
\left \vert X_k \right \vert_{g_k(\tau)} \leq C_0 \left ( \sup_{
D( r_1) \times \left[ 0,T \right] }
\left \vert X_\infty \right \vert_{g_\infty(\tau)}  +1 \right ) \notag\\
&  \sup_{\Phi_k \left(D( r_1) \right) \times \left[ 0,T \right] }
\max \{ -Q_k , 0 \}  \leq C_0 \left ( \sup_{ D( r_1)
 \times \left[ 0,T \right] }
\max \{ -Q_\infty , 0 \}   +1 \right ) \notag.
\end{align}
Note that there is a constant $K \geq 0$
such that
\begin{equation}
\operatorname{Rc}\left( \tilde{g}\right) \geq -K \,\,\,\, \text{ in } D \left(
r_1\right). \label{eq K Ricci bdd}
\end{equation}

To apply Theorem \ref{thm mvi} to $u=H_k$ we choose $\Omega$ in Theorem
\ref{thm mvi} to be $D(2r_0)$ where
\[
r_{0} \in (0, \min \left\{ \frac{\sqrt{\tau _{1}}}{2},1, \frac{r_1}{2}\right\}) .
\]
we take $g(\tau) =g_k(\tau)$, $X=X_k$, and $Q=Q_k$.
Now we verify the assumption of Theorem \ref{thm mvi}.
For any $\left( x_{\ast },\tau _{\ast }\right) \in D \times
\left[ \tau _{1},\tau _{2}\right] $,
by the choice of $r_1$ and $r_0$ we conclude that the closed ball
$\overline{B_{\tilde{g}}(x_*, 2r_0) }$ is compact in $D(r_1)$
and that $\tau_* \geq (2r_0)^2$.
Also for any $x_{\ast } \in D $ we have
\begin{equation*}
\operatorname{Vol}_{\tilde{g}}\left( B_{\tilde{g}}\left( x_{\ast },r_{0}\right)
\right) \geq c_{1}\operatorname{Vol}_{\tilde{g}}\left( B_{\tilde{g}}\left(
x_{0\infty },r_{0}\right) \right)
\end{equation*}
where $c_{1}$ is a constant depending on $n$, $K$, $r_{0}$ and
$\operatorname{diam}_{\tilde{g}}\left( D \right) $,
this follows from a standard argument using Bishop-Gromov volume comparison theorem.

By the local mean value property (\ref{eq mvi})
and bounds in (\ref{eq get k rid of})
we have the following. When $k \geq k_0$, for any $\left(
x_{\ast },\tau _{\ast }\right) \in D \times \left[
\tau _{1},\tau _{2}\right] $ and for any $\tilde{x} \in B_{\tilde{g}} (x_\ast, r_{0})$
and $\tilde{\tau} \in [\tau_\ast -r_0^2,\tau_\ast]$
\begin{eqnarray*}
 0 &\leq & H_{k}\left(\Phi_k(\tilde{ x} ), \tilde{\tau} \right)
\leq \sup_{B_{ \tilde{g}}\left( x_{\ast }, r_0 \right ) \times [\tau_{\ast }
-(r_{0})^2, \tau_\ast ] } H_{k}\left(\Phi_k( x) ,\tau \right) \\
& \leq & C_4 \int_{B_{
\tilde{g}}\left( x_{\ast }, 2r_0 \right ) \times [\tau_{\ast }
-(2r_{0})^2, \tau_\ast ] }
H_{k}\left( \Phi_k(y), \tau \right) d\mu _{\tilde{g}}\left(
y \right) d\tau  \\
& \leq & C_4 \int_{B_{
\tilde{g}}\left( x_{\ast }, 2r_0 \right ) \times [\tau_{\ast }
-(2r_{0})^2, \tau_\ast ]  } H_{k}\left( \Phi_k(y),\tau \right) \left(
C_0^{n}d\mu _{\Phi_k^* g_{k}\left( \tau
\right) }\left( y \right) \right) d\tau  \\
& \leq & C_4 C_0^{n}\int_{\tau _{\ast }-\left( 2r_{0}\right)
^{2}}^{\tau _{\ast }}\int_{\Phi_k ( B_{
\tilde{g}}\left( x_{\ast }, 2r_0 \right ) )}
H_{k}\left( z,\tau \right) d\mu _{g_{k}\left( \tau \right)
}\left( z \right) d\tau,
\end{eqnarray*}
where
\[
C_4 \doteqdot  \frac{C_{1}e^{C_{2}T+C_{3}\sqrt{K}r_{0}}}{r_{0}^{2}
\operatorname{Vol}_{\tilde{g}}\left( B_{\tilde{g}}\left( x_{\ast },
r_{0}\right) \right) }
\leq  \frac{C_{1}e^{C_{2}T+C_{3}\sqrt{K}r_{0}}}{r_{0}^{2} c_1
\operatorname{Vol}_{\tilde{g}}\left( B_{\tilde{g}}\left( x_{0\infty},
r_{0}\right) \right) }.
 \]
Here constant $C_{1}$ depends only on $n$, $T$, $C_0$,
and $\sup_{D(r_1) \times \left[ 0,T \right] }\left\vert \left( \mathcal{R}_{
\infty }\right) _{ij} \right\vert _{g_{\infty }\left( \tau \right) }$.
Constant $C_{2}$ depends only on $C_0$,
$\sup_{D(r_1) \times \left[ 0,T\right] }\max \{-Q
_{\infty }, 0 \}$, $\sup_{D(r_1) \times \left[ 0,T\right] }$ $
\left \vert X_{\infty } \right\vert_{g_\infty (\tau)}$,
and $\sup_{D(r_1) \times \left[ 0,T\right] }$ $\left\vert \mathcal{R}
_{\infty } \right\vert $. Constant $C_{3}$ depends only on
$n$.

In summary so far we have proved that when $k \geq k_0$, for any $\left(
x,\tau \right) \in D(r_0) \times \left[
\tau _{1}-r_0^2, \tau _{2}\right] $
\begin{equation}
\left\vert H_{k}\left( \Phi_k( x),\tau \right) \right\vert \leq
 C_5 \int_{\tau _1-\left( 2r_{0}\right)
^{2}}^{\tau_2 }\int_{\Phi_k ( B_{
\tilde{g}}\left( x_{\ast }, 2r_0 \right ) )}
H_{k}\left( z,\tau \right) d\mu _{g_{k}\left( \tau \right)
}\left( z \right) d\tau, \label{eq what is what use 01}
\end{equation}
where $C_5$ is a constant independent of $k$ and $x_\ast$ is a point in $D$
such that $x \in B_{\tilde{g}}(x_\ast, r_0)$.
It follows from the assumption (\ref{eq Hk int assump}) that
\begin{equation}
\sup_{D(r_0) \times \left[ \tau _{1}-r_0^2, \tau _{2}\right]}
\left\vert H_{k}\left( \Phi_k( x),\tau \right) \right\vert \leq C_6
\label{eq Hk c0 bound}
\end{equation}
where $C_6$ is a constant independent of $k$.

%%%%%%%%%%%%%%%%%%%%%%%%%%%%%%%%%%%%%%%%%%%%%%%%%%%%%%%%%%%%%%%%%%%%%%%%%%%%
%%%%%%%%%%%%%%%%%%%%%%%%%%%%%%%%%%%%%%%%%%%%%%%%%%%%%%%%%%%%%%%%%%%%%%%%%%%%
%%%%%%%%%%%%%%% \newpage
\vskip .5cm
Next we want to get local high derivative estimates of $H_k ( \Phi_k(\cdot), \cdot)$.
Below we assume $k \geq k_0$.
From the compactness of $\overline{D(r_1)}$ there is a $K_1 >0$ such that
$\left \vert \operatorname{Rm}_{\tilde{g}} (x) \right \vert_{\tilde{g}}
\leq K_1$ for all $x \in D(r_1)$.
Let $\hat{r}_0\doteqdot \min \left\{ r_{0},\frac{\pi }{4\sqrt{K_{1}}}\right\}$.
Fix an arbitrary $ x_{\ast } \in D $,
let $\exp _{x_{\ast }}:B\left( 0,\hat{r}_0 \right) \rightarrow B_{\tilde{g}
}\left( x_{\ast },\hat{r}_0 \right) \subset M_\infty$ be the exponential
map of metric $\tilde{g}$ and let $\vec{x}=\left( x_{i}\right)$ be an associated
normal coordinates on $B\left( 0,\tilde{r}_0 \right)$.
We will consider the pull back functions
\begin{equation}
\hat{H}_k(\vec{x},\tau ) \doteqdot H_{k}\left( \Phi_k(\exp_{x_\ast} \vec{x} ,
\tau \right)  \label{eq uk def}
\end{equation}
defined on $B\left( 0,\hat{r}_0 \right) \times (0,T]$.

Let $\hat{X}_k( \tau)$ be the vector fields on $B\left( 0,\hat{r}_0 \right)$
such that the push-forward vector field $(\Phi_k \circ \exp_{x_\ast} )_*  \hat{X}_k(\tau)= X_k( \tau)$.
Since $H_k$ satisfies (\ref{eq seq fund sol}),
$\hat{H}_k$ satisfies the following
\begin{equation}
\frac{\partial \hat{H}_k}{\partial\tau}- \Delta_{(\Phi_k \circ \exp_{x_\ast} )^* g_k
\left(  \tau\right)  }\hat{H}_k + \nabla_{\tilde{X}_k} \hat{H}_k + Q_k(\Phi_k(\exp_{x_\ast}
\vec{x} ),\tau ) \hat{H}_k  = 0 \label{eq uks equ compsed}
\end{equation}
on $B\left( 0,\hat{r}_0 \right)  \times (0,T]$.

Because of the $C^\infty$ Cheeger-Gromov convergence of
$( \Omega_k, g_k(\tau),(x_{0k},0) ) $
to $ ( M_\infty, g_{\infty }(\tau), (x_{0\infty },0) )$
for $\tau \in \left[ 0,T \right]$,
the structure coefficient functions of the parabolic equation
(\ref{eq uks equ compsed})
converge in $C^{\infty}$-norm on $B\left( 0,\hat{r}_0 \right) \times [0,T]$
to the corresponding coefficients of the following parabolic equation
\begin{equation}
\frac{\partial \hat{H}_\infty }{\partial\tau}- \Delta_{\exp_{x_\ast}^* g_\infty
\left(  \tau\right)  }\hat{H}_\infty+ \nabla_{\hat{X}_\infty}
\hat{H}_\infty + Q_\infty (\exp_{x_\ast}
\vec{x} , \tau ) \hat{H}_\infty  = 0, \label{eq uinfty equ compsed}
\end{equation}
where  $\hat{X}_\infty (\cdot,\tau)$ is the vector fields
on $B\left( 0,\hat{r}_0 \right)$
such that $(\exp_{x_\ast} )_*  \hat{X}_\infty (\cdot,\tau)= X_\infty(\cdot,\tau)$.
In particular for any fixed $l \in \mathbb{N}$ and $\alpha \in (0,1)$
the parabolic H\"{o}lder norms
\[
\left \Vert  \cdot \right \Vert_{C^{2l,\alpha; l,\alpha/2} (
B\left( 0,\hat{r}_0 \right) \times [0,T]) }
\]
of these structure coefficient functions are bounded by a constant
independent of $k$.
Note that from (\ref{eq Hk c0 bound}) $\hat{H}_{k}$ is
uniform bounded on $B \left( 0,\hat{r}_0 \right) \times
\left[ \tau_{1}-r_0^2, \tau_{2} \right] $.

We can apply the interior Schauder estimates for linear
parabolic equation (see
Theorem 4.9 and Exercise 4.5 in \cite{Lie}, for example) to conclude that for any integer $l \geq 0$ and  $\alpha \in (0, 1) $ we have
\begin{equation}
\left \Vert \hat{H}_{k} \right \Vert_{C^{2l+2, \alpha ;l+1,\alpha/2} \left( \bar{B}
\left( 0, \hat{r}_0/2 \right) \times \left[ \tau_1,\tau_2 \right ] \right )
} \leq C_{l,\alpha} \label{eq holder est uk}
\end{equation}
where $C_{l,\alpha}$ is a constant independent of $k$.

Finally we can show the existence of $H_\infty$ in Theorem \ref{thm de main conv}.
By estimates (\ref{eq holder est uk}) and Arzela and Ascoli theorem there is a subsequence (still indexed by $k$) such that $\hat{H}_{k}\rightarrow \hat{H}
_{\infty }$ in $C^\infty$-norm on $ \bar{B}\left( 0,\hat{r}_0/2 \right)
\times \left[ \tau_{1},\tau _{2} \right]$ and
$\hat{H}_{\infty }$ is a solution of (\ref{eq uinfty equ compsed}).
Note that $\tilde{H}_k(x,\tau) = H_{k}( \Phi_k(x), \tau)$ satisfies the
following equation
\begin{equation}
\frac{\partial \tilde{H}_k}{\partial\tau}- \Delta_{\Phi_k^* g_k
\left(  \tau\right)  }\tilde{H}_k + \nabla_{\tilde{X}_k} \tilde{H}_k + Q_k(\Phi_k(x),\tau )
\tilde{H}_k  = 0 \label{eq tilde H k}
\end{equation}
where $\tilde{X}_k(\tau)$ is a family of vector field such that the push-forward
vector filed $(\Phi_k)_* \tilde{X}_k ( \tau) =X_k( \tau)$.
Define function $H_\infty$ on $\exp_{x_\ast} \bar{B}\left( 0,\hat{r}_0/2 \right)
\times [\tau_1, \tau_2]$ by $H_{\infty }(\exp_{x_\ast}(\vec{x}),\tau) \doteqdot \hat{H}_\infty (\vec{x} ,\tau) $.
Though $\exp_{x_\ast}$ is not necessarily a diffeomorphism,
it is clear that $H_\infty$ is well defined and that
$\tilde{H}_k$ converges to $H_\infty$ in $C^\infty$-norm.
Furthermore $H_\infty$ is a solution of (\ref{eq Hinfty sat})
on  $\exp_{x_\ast} \bar{B}\left( 0,\hat{r}_0/2 \right) \times [\tau_1, \tau_2] $.

Note that our definition of $ \hat{H}_k(\vec{x},\tau)$ in (\ref{eq uk def})
depends on $x_\ast$.
Since $x_{\ast} $ is an arbitrary point in $D$,
we conclude by a diagonalization argument that there is a
subsequence of $\tilde{H}_k$
which converges in $C^\infty$-norm to some function $H_{\infty }$
defined on $D \times \left[ \tau_1,\tau_{2} \right]$.
$H_{\infty }$ is a solution of (\ref{eq Hinfty sat})
on $D\times \left[ \tau_1,\tau_{2} \right]$.

Let $\tau_2 =T$ and let $\tau_{1i}$ be a sequence approaching to $0^+$.
Let $D_i$ be a sequence of compact domains
which exhaust $\operatorname{Int}(M_\infty)$.
We can apply the above construction of $H_\infty$ on
$D\times \left[ \tau_1,\tau_{2} \right]$ repeatedly to each
$D_i\times \left[ \tau_{1i}, T \right]$.
Using a diagonalization argument we can find a subsequence such that
$ \tilde{H}_{k}$ converges in $C^\infty$-norm on $D_i\times
\left[ \tau_{1i}, T \right]$ for each $i$.
This implies that $ \tilde{H}_{k}$ converges on
$\operatorname{Int}(M_{\infty }) \times ( 0,T ]$ to a function $H_\infty$,
the convergence is uniform in $C^\infty$-norm on any compact subset
of $\operatorname{Int}(M_{\infty }) \times ( 0,T ]$.

In short we have defined a function $H_{\infty } \geq 0$
 which satisfies (\ref{eq Hinfty sat}) on $\operatorname{Int}(M_{\infty })
  \times ( 0,T ]$.

%%%%%%%%%%%%%%%%%%%%%%%%%%%%%%%%%%%%%%%%%%%%%%%%%%%%%%%%%%%%%%%%%%%%%%%%%%%%
%%%%%%%%%%%%%%%%%%%%%%%%%%%%%%%%%%%%%%%%%%%%%%%%%%%%%%%%%%%%%%%%%%%%%%%%%%%%
%%%%%%%%%\newpage
\vskip .5cm
\noindent \textbf{2.5 Step 2 Proof of (\ref{eq delta thm main 1})}.
We will show $\lim_{\tau \rightarrow 0^+}
H_\infty(x,\tau) = \delta_{x_{0\infty}}$.
Let $F: M_\infty \rightarrow \mathbb{R}$ be a nonnegative
 $C^2$ function such that support
$\operatorname{supp}(F)$ is a compact subset of $\operatorname{Int}(M_\infty)$.
We compute using (\ref{eq tilde H k})

\begin{align}
&  \frac{d}{d\tau}\int_{M_{\infty}}\tilde{H}_{k}Fd\mu_{\Phi_{k}^{\ast}%
g_{k}(\tau)} \notag \\
= &  \int_{M_{\infty}}\left(  \frac{\partial\tilde{H}_{k}}{\partial\tau
}+\mathcal{R}_{k}(\Phi_{k}(x),\tau)\tilde{H}_{k}\right)  Fd\mu_{\Phi_{k}%
^{\ast}g_{k}(\tau)} \nonumber \\
= &  \int_{M_{\infty}}\left(
\begin{array}
[c]{c}%
\Delta_{\Phi_{k}^{\ast}g_{k}(\tau)}\tilde{H}_{k} - \nabla_{\tilde{X}_{k}}
\tilde{H}_{k}\\
+\left(  \mathcal{R}_{k}(\Phi_{k}(x),\tau)-{Q}_{k}(\Phi_{k}(x),\tau)\right)
\tilde{H}_{k}%
\end{array}
\right)  Fd\mu_{\Phi_{k}^{\ast}g_{k}(\tau)}\nonumber\\
= &  \int_{M_{\infty}}\left(
\begin{array}
[c]{c}%
\left ( \Delta_{\Phi_{k}^{\ast}g_{k}(\tau)}F \right ) \tilde{H}_{k} +
\operatorname{div} ( F \tilde{X}_k ) \tilde{H}_{k} \\
+\left(  \mathcal{R}_{k}(\Phi_{k}(x),\tau)-{Q}_{k}(\Phi_{k}(x),\tau)\right)
\tilde{H}_{k}F
\end{array}
\right)  d\mu_{\Phi_{k}^{\ast}g_{k}(\tau)}, \label{eq uk F integ deri bdd}%
\end{align}
where in the last equality we have used the divergence theorem
\begin{align*}
 \int_{M_\infty} \left ( -\nabla_{\tilde{X}_k} \tilde{H}_k \right ) F
d\mu_{\Phi_{k}^{\ast}g_{k}(\tau)} 
= & \int_{M_\infty}  - \left \langle F \tilde{X}_k, \nabla\tilde{H}_k \right \rangle_{\Phi_{k}^{\ast}g_{k}(\tau)} d\mu_{\Phi_{k}^{\ast}g_{k}(\tau)} \\
= & \int_{M_\infty} \operatorname{div}_{\Phi_{k}^{\ast}g_{k}(\tau)}
( F \tilde{X}_k ) \tilde{H}_{k} d\mu_{\Phi_{k}^{\ast}g_{k}(\tau)}.
\end{align*}

Because of the following uniform convergence on
$\operatorname{supp}(F) \times [0,T]$
\begin{align*}
& \Delta_{\Phi_k^* g_k(\tau)} F
\rightarrow \Delta_{ g_\infty (\tau)} F \\
&  \operatorname{div}_{\Phi_{k}^{\ast}g_{k}(\tau)} ( F \tilde{X}_k ) \rightarrow
\operatorname{div}_{g_\infty(\tau)} ( F X_\infty ) \\
& \mathcal{R}_k(\Phi_k(x),\tau) \rightarrow \mathcal{R}_\infty(x,\tau) \\
& {Q}_k(\Phi_k(x),\tau) \rightarrow {Q}_\infty(x,\tau) ,
\end{align*}
there is a constant $C_7$ independent of $k$ such that
\begin{align*}
& \sup_{ \operatorname{supp}(F) \times [0,T]} \vert
\Delta_{\Phi_k^* g_k(\tau)} F \vert \leq C_7  \\
& \sup_{ \operatorname{supp}(F) \times [0,T]} \vert \operatorname{div}_{\Phi_{k}^{\ast}g_{k}(\tau)} ( F \tilde{X}_k )
\vert \leq C_7 \\
& \sup_{ \operatorname{supp}(F) \times [0,T]} \vert
\mathcal{R}_k(\Phi_k(x),\tau)  \vert \leq C_7 \\
& \sup_{ \operatorname{supp}(F) \times [0,T]} \vert
{Q}_k(\Phi_k(x),\tau) \vert \leq C_7.
\end{align*}
Here and below we assume $k$ are large enough.
Hence it follows from (\ref{eq uk F integ deri bdd}) and (\ref{eq Hk int assump})
\begin{equation}
\left \vert  \frac{d}{d \tau} \int_{M_\infty} \tilde{H}_k F
d\mu_{\Phi_k^* g_k(\tau)} \right \vert \leq 2 C_7  \int_{M_\infty}
\tilde{H}_k F d \mu_{\Phi_k^* g_k(\tau)} + 2 C_8 \label{eq uk F ode bdd}
\end{equation}
for $\tau \in ( 0,T]$ where $C_8 \doteqdot C_7 C_*$.

Let
\[
U_k(\tau) \doteqdot \int_{M_\infty} \tilde{H}_k(x,\tau) F(x)
d \mu_{\Phi_k^* g_k(\tau)}(x).
\]
By a simple integration of (\ref{eq uk F ode bdd}) we get
\begin{equation}
 U_k(0) e^{-2C_7 \tau}  - \frac{C_8}{C_7} \left ( 1 - e^{-2C_7 \tau}  \right )
  \leq  U_k(\tau) \leq
U_k(0) e^{2C_7 \tau}   +\frac{C_8}{C_7} \left ( e^{2 C_7 \tau}
 -1 \right ) . \label{eq uk F upp lower bdd}
\end{equation}
By the definition of $\tilde{H}_k$ we have
\[
U_k(0) = \int_{M_{\infty }} H_{k}(y,0) F(\Phi
_{k}(y))d\mu_{g_{k}\left( \tau \right) }  = F(x_{0\infty}) .
\]
By the convergence of $\tilde{H}_k$ proved in \S 2.4 and the compactness of
 $ \operatorname{supp}(F)$ we have that for any $\tau >0$
\begin{align*}
\lim_{k \rightarrow \infty} U_k(\tau)
& =  \int_{M_{\infty }} \lim_{k \rightarrow \infty} \left ( \tilde{H}
_{k}(x,\tau)F(x) d\mu _{\Phi _{k}^{\ast }g_{k}\left(
\tau \right) }(x) \right )   \\  %%%\notag
& =  \int_{M_{\infty }}H_{\infty}(x,\tau) F(x) d\mu_{g_\infty \left(
\tau \right) } (x). %%%  \label{eq u k tau limit}
\end{align*}
By taking $k \rightarrow \infty$ limit of (\ref{eq uk F upp lower bdd}) we get
\begin{align*}
e^{-2C_7 \tau} F(x_{0\infty})- \frac{C_8}{C_7} \left ( 1 - e^{-2C_7 \tau}  \right )
&  \leq  \int_{M_{\infty }}H_{\infty} F d\mu_{g_\infty \left(
\tau \right) } \\
&  \leq  e^{2C_7 \tau} F(x_{0\infty})  +\frac{C_8}{C_7} \left ( e^{2
C_7 \tau} -1 \right )
\end{align*}
for $\tau \in ( 0,T]$.
Hence
\begin{equation}
\lim_{\tau \rightarrow 0^+} \int_{M_{\infty }}H_{\infty} F d\mu_{g_\infty \left(
\tau \right) } = F(x_{0\infty}) \label{eq H infty delta prop}
\end{equation}
for any nonnegative $C^2$ function with compact support in $\operatorname{Int}(M_\infty)$. This implies $\lim_{\tau \rightarrow 0^+}
H_\infty(x,\tau) = \delta_{x_{0\infty}}$.

%%%%%%%%%%%%%%%%%%%%%%%%%%%%%%%%%%%%%%%%%%%%%%%%%%%%%%%%%%%%%%%%%%%%%%%%%%%%
%%%%%%%%%%%%%%%%%%%%%%%%%%%%%%%%%%%%%%%%%%%%%%%%%%%%%%%%%%%%%%%%%%%%%%%%%%%%
%%%%%%%%%\newpage
\vskip .2cm
\noindent \textbf{2.6 Step 3 Finishing the proof of Theorem \ref{thm de main conv}}.
By (\ref{eq H infty delta prop}) we conclude that there is a $T_1 \in (0,T)$
 such that for any $\tau \in (0,T_1]$
there is a $x_\tau \in B_{\tilde{g} } \left( x_{0\infty}, \hat{r}_0 \right) $
such that $H_{\infty }\left( x_\tau, \tau \right) >0$.
Since $H_{\infty }\geq 0$ and $H_\infty$ satisfies the equation
(\ref{eq Hinfty sat}) on $M_\infty \times(0,T]$,
it follows from the strong maximum principle
that $H_\infty(x,\tau) >0$ on $M_\infty \times(0,T]$.
Hence $\ln \tilde{H}_{k}\left(x,\tau \right) \rightarrow \ln
H_{\infty }\left( x,\tau \right) $ uniformly in
$C^{\infty }$-norm on any compact subset of $\operatorname{Int}(M_{\infty })
\times (0,T]$,
and the claimed convergence of $\tilde{f}_k$ to $f_\infty$
in Theorem \ref{thm de main conv} follows easily.
Now Theorem \ref{thm de main conv} is proved.

%%%%%%%%%%%%%%%%%%%%%%%%%%%%%%%%%%%%%%%%%%%%%%%%%%%%%%%
%%%%%%%%%%%%%%%%%%%%%%%%%%%%%%%%%%%%%%%%%%%%%%%%%%%%%%%
\section{The upper bound of integrals of fundamental solutions}

In this section we discuss how to bound the integrals of fundamental solutions,
this is related to the assumption (\ref{eq Hk int assump}).
Here we use a special case of the setup described at the beginning of
\S 2.1, we assume $\Omega =M$. We also adopt the notations used at the
beginning of \S 2.1, in particular, $H$ is a fundamental solution of
(\ref{eq fund sol}) on $M \times [0,T]$ centered at $x_0 \in \operatorname{Int}(M)$.
It is clear
\begin{equation}
\lim_{\tau \rightarrow 0^+} \int_M H(x,\tau) d\mu_{g(\tau)} (x) =1.
\label{eq H int time 01}
\end{equation}

The results in this section and the next section are well-known,
we include them here because of their relation with Theorem \ref{thm de main conv}.
Below we divide our discussion into three cases.

\vskip .2cm
\noindent \textbf{Case 1 $M$ is a closed manifold}. Let $C_0$ be a constant
such that
\[
\sup_{M \times [0,T]} (\operatorname{div}_{g(\tau)} X +
\mathcal{R}-Q ) \leq C_0.
\]
We compute
\begin{align}
\frac{d}{d\tau} \int_M H d\mu_{g(\tau)} & = \int_M \left ( \frac{\partial
H}{\partial \tau } + \mathcal{R} H \right ) d\mu_{g(\tau)} \notag  \\
& = \int_M \left ( \Delta H - \nabla_X H + (\mathcal{R}-Q )
H \right ) d\mu_{g(\tau)}  \notag \\
& =\int_M  (\operatorname{div}_{g(\tau)} X+ \mathcal{R}-Q ) H
 d\mu_{g(\tau)}  \label{eq der heat kernel int}  \\
& \leq C_0 \int_M H  d\mu_{g(\tau)}.  \notag
\end{align}
Hence it follows from (\ref{eq H int time 01}) that
\begin{equation*}
\int_M H d\mu_{g(\tau)} \leq e^{ C_0 \tau}
\end{equation*}
for $\tau \in [0,T]$.

In the special case when $X=0$ and $Q =\mathcal{R}$ it follows from
(\ref{eq der heat kernel int}) that
\[
\int_M H d\mu_{g(\tau)}=1
\]
for $\tau \in [0,T]$.

%%%%%%%%%%%%%%%%%%%%%%%%%
\vskip .2cm
\noindent \textbf{Case 2 $M$ is a compact manifold with nonempty boundary}.
Let $C_0$ be a constant
such that
\[
\sup_{M \times [0,T]} (\operatorname{div}_{g(\tau)} X +\mathcal{R}-Q )
 \leq C_0.
\]
Let $\nu=\nu(x,\tau)$ be the outward unit normal direction on $\partial M$ with
respect to the metric $g(\tau)$.
We assume $\langle X, \nu \rangle_{g(\tau)} \geq 0$ on the boundary $\partial M
\times [0,T]$.
We compute
\begin{align}
& \frac{d}{d\tau} \int_M H d\mu_{g(\tau)} \notag \\
=& \int_M \left ( \Delta H
- \nabla_X H + ( \mathcal{R}-Q ) H \right ) d\mu_{g(\tau)} \notag  \\
 =& \int_{\partial M} \left ( \frac{\partial H}{\partial \nu} - \langle X, \nu
\rangle_{g(\tau)} H  \right ) d \mu_{g(\tau)}
 + \int_{ M}  (\operatorname{div}_{g(\tau)} X + \mathcal{R}-Q ) H  d\mu_{g(\tau)}
 \notag \\
 \leq & \int_{\partial M}  \frac{\partial H}{\partial \nu}  d\mu_{g(\tau)}
 + \int_{ M}  (\operatorname{div}_{g(\tau)} X + \mathcal{R}-Q ) H  d\mu_{g(\tau)} .
 \label{eq der heat kernel int 2}
\end{align}

If $H$ satisfies the Dirichlet boundary condition $\left . H \right
\vert_{\partial M \times (0,T]} =0$, then $ \frac{\partial H}{\partial \nu}
\leq 0$ on $\partial M \times (0,T]$.
It follows from (\ref{eq der heat kernel int 2}) that
\[
\int_M H d\mu_{g(\tau)} \leq e^{ C_0 \tau}
\]
for $\tau \in (0,T]$. In the special case when $X=0$ and
 $Q =\mathcal{R}$ it follows  from (\ref{eq der heat kernel int 2}) that
\[
\int_M H d\mu_{g(\tau)} \leq 1 \,\,\, \text{ for any } \tau \in (0,T].
\]

If $H$ satisfies the Neumann boundary condition $\left . \frac{\partial
H}{\partial \nu} \right \vert_{\partial M \times (0,T]} =0$,
then  (\ref{eq der heat kernel int 2}) becomes
\[
\frac{d}{d\tau} \int_M H d\mu_{g(\tau)} \leq
\int_{\partial M}  (\operatorname{div}_{g(\tau)} X
+\mathcal{R}-Q ) H  d\mu_{g(\tau)}.
\]
Hence
\[
\int_M H d\mu_{g(\tau)} \leq e^{ C_0 \tau} \,\,\, \text{ for any } \tau \in ( 0,T].
\]
In the special case when $X=0$ and $Q =\mathcal{R}$
it follows that
\[
\int_M H d\mu_{g(\tau)} = 1 \,\,\, \text{ for any } \tau \in (0,T].
\]

%%%%%%%%%%%%%%%%%%%%%%%%%%%%%%%%%%%%%%%%%%%%%%%%%%%%%%%%%%%%%%%%%%%%%%%%%%%%%%
\vskip .2cm
\noindent \textbf{Case 3 $(M,g(\tau)), \, \tau \in [0,T]$, are complete
noncompact manifolds}. Because of the potential non-uniqueness of the fundamental
solutions, it is desirable to consider the integrals of the minimal fundamental solutions. We refer the reader to the literature (see, for example, Lemma 5.1
in \cite{CTY}, or Lemma 26.14 and Corollary 26.15 in \cite{CCIII}).

%%%%%%%%%%%%%%%%%%%%%%%%%%%%%%%%%%%%%%%%%%%%%%%%%%%%%%%%%%%%%%%%%%%%%%%%%%%%
%%%%%%%%%%%%%%%%%%%%%%%%%%%%%%%%%%%%%%%%%%%%%%%%%%%%%%%%%%%%%%%%%%%%%%%%%%%%
%%%%% \newpage
\section{Uniqueness of fundamental solutions}

In this section we consider the uniqueness of fundamental solutions.
Note that it is pointed out by Hsu in \cite{Hs} that when the limit fundamental
solution in Theorem \ref{thm de main conv} is unique, then
the sub-convergence in the theorem can be improved to be convergence
for the whole sequence.
Here we use the same setup as the one used \S 3,
in particular, $H$ is a fundamental solution of
(\ref{eq fund sol}) on $M^n \times [0,T]$ centered at $x_0
\in \operatorname{Int}(M)$.

\begin{proposition}
(i) When $M$ is a closed manifold, the fundamental solution is unique.

\noindent (ii) When $M$ is a compact manifold with nonempty boundary,
the fundamental solution with any Dirichlet boundary condition is unique.

\noindent (iii) When $M$ is a compact manifold with nonempty boundary,
the fundamental solution with any Neumann boundary condition is unique.

\noindent (iv) When $(M,g(\tau)), \, \tau \in [0,T]$, are complete
and noncompact manifold,
the fundamental solution satisfying the assumption (A1), (A2), and
(A3) below is unique.
\end{proposition}

\noindent \textbf{Remark}
The proof is based on an idea of Brett Kotschwar (see Footnote 14 on
 p.345 of \cite{CCII}).

\vskip .1cm
\begin{proof}
(i) Let $H_1$ and $H_2$ be two fundamental solutions centered at $x_0$.
Define $F(x,\tau) \doteqdot H_1(x,\tau) - H_2(x,\tau)$.
Then $F$ satisfies
\begin{equation}
\frac{\partial F}{\partial\tau}
- \Delta_{g\left(  \tau\right)  }F + \nabla_X F + Q F =0.
\label{eq F difference unique}
\end{equation}
Let $\varphi: M \rightarrow \mathbb{R}$ be an arbitrary $C^{2}$ function.
Fix a $\bar{\tau} \in (0,T]$.
Let $\Phi: M \times\left[ 0,\bar{\tau}\right]  \rightarrow \mathbb{R}$
be the solution to the initial value problem
\begin{align}
& \left(  \frac{\partial}{\partial\tau}+\Delta_{g\left(  \tau\right)  }
+ \nabla_X +(\operatorname{div}_{g(\tau)} X +\mathcal{R}-Q) \right) \Phi   =0,
\notag \\
& \Phi\left(  \cdot,\bar{\tau}\right)  =\varphi.  \label{eq Phi sol backward eq test}
\end{align}
The solution always exists.

For $\varepsilon\in\left( 0,\bar{\tau}\right)$ we have
\begin{align*}
 0    =& \int_{\varepsilon}^{\bar{\tau}}\int_{ M }\Phi
\left(  \frac{\partial}{\partial\tau}-\Delta_{g\left(
\tau\right)  } + \nabla_X +Q \right)  F  d\mu_{g(\tau) } d\tau  \\
= & \int_{\varepsilon}^{\bar{\tau}}\left(  \frac{d}{d\tau}\int_{M
}\Phi F d\mu_{g(\tau) }\right)  d\tau \\
& -\int_{\varepsilon}^{\bar{\tau}}
\int_{ M}F\left(  \frac{\partial}{\partial\tau}+\Delta_{g\left(
\tau\right) } + \nabla_X + (\operatorname{div}_{g(\tau)} X
+\mathcal{R}-Q)  \right)  \Phi d\mu_{g(\tau)}d\tau\\
 = & \int_{M }\Phi Fd\mu_{g(\bar{\tau})}-\int_{M}\Phi
Fd\mu_{g(\varepsilon)}.
\end{align*}
Since
\begin{align*}
 \lim_{\varepsilon \rightarrow 0} \int_{M}\Phi Fd\mu_{g(\varepsilon)}
& = \lim_{\varepsilon \rightarrow 0} \int_{M}\Phi H_1 d\mu_{g(\varepsilon)}
- \lim_{\varepsilon \rightarrow 0} \int_{M}\Phi H_2
d\mu_{g(\varepsilon)} \\
& =\Phi(x_0,0) -\Phi (x_0,0) = 0,
\end{align*}
we have proved
\[
\int_{M }\varphi F d\mu_{g(\bar{\tau})} = \int_{M }\Phi F d\mu_{g(\bar{\tau})} = 0.
\]
Since both $\varphi$ and $\bar{\tau}$ are arbitrary,
we concluded $F=0$ and hence $H_1 =H_2$.

%%%%%%%%%%%%%%%%%%%%%%%%%%%%%%%%%%%%%%
\vskip .1cm
(ii) Let $\psi:M \times [0,T] \rightarrow \mathbb{R}$ be a continuous function.
Let $H_1$ and $H_2$ be two fundamental solution of (\ref{eq fund sol})
on $M \times (0,T]$ centered at $x_0$ which satisfy the following
Dirichlet boundary condition
\[
\left . H_1 \right \vert_{\partial M \times (0,T]} =
\left . H_2 \right \vert_{\partial M \times (0,T]}
=\left . \psi \right \vert_{\partial M \times (0,T]}.
\]
Define $F(x,\tau) \doteqdot H_1(x,\tau) - H_2 (x,\tau)$.
Then $F$ satisfies (\ref{eq F difference unique}) and the boundary condition
$\left . F \right \vert_{\partial M \times ( 0,T]} = 0$.
Let $\varphi: M \rightarrow \mathbb{R}$ be an arbitrary $C^{2}$ function
which vanishes on $\partial M$.
Fix a $\bar{\tau} \in (0,T]$.
Let $\Phi: M \times\left[ 0,\bar{\tau}\right]  \rightarrow \mathbb{R}$
be the solution to the initial-boundary value problem
(\ref{eq Phi sol backward eq test}) with Dirichlet boundary condition
$\left . \Phi \right \vert_{ \partial M \times [0,T]} =0$.
The solution always exists.

By the divergence theorem we have
\begin{align}
& \int_M \Phi \Delta_{g(\tau)} F   d \mu_{g(\tau)} = &
  \int_M  F  \Delta_{g(\tau)}\Phi  d \mu_{g(\tau)}
  + \int_{\partial M} \left ( \frac{\partial F}{\partial \nu} \Phi
  - F \frac{\partial \Phi }{\partial \nu} \right ) d\sigma_{g(\tau)} \notag \\
  &  \int_M  \Phi \nabla_X F  d \mu_{g(\tau)} = & \int_M \left (- F \nabla_X  \Phi
  -  F  \Phi \operatorname{div}_{g(\tau)} X   \right )d \mu_{g(\tau)} \notag \\
  & & + \int_{\partial M} \langle X, \nu \rangle  F \Phi
d\sigma_{g(\tau)}  \label{eq int by parts F}
 \end{align}
where $d\sigma_{g(\tau)}$ is the volume form on $\partial M$ defined by the metric
$\left . g(\tau) \right \vert_{\partial M}$.
Because $F= \Phi =0$ on $\partial M \times (0,T]$,
it follows from similar calculations and arguments as in (i) that
for $\varepsilon\in\left( 0,\bar{\tau}\right)$
\[
\int_{M }\Phi Fd\mu_{g(\bar{\tau})}= \int_{M}\Phi Fd\mu_{g(\varepsilon)}
\]
and hence $F=0$. This proves $H_1 =H_2$.

%%%%%%%%%%%%%%%%%%%
%%%%%%%%%%%%%%%%%%%
\vskip .1cm
(iii) Let $\psi:M \times [0,T] \rightarrow \mathbb{R}$ be a continuous function.
Let $H_1$ and $H_2$ be two fundamental solution of (\ref{eq fund sol})
on $M \times (0,T]$ centered at $x_0 \in \operatorname{Int}(M)$
which satisfy the following
Neumann boundary condition
\[
\left . \frac{\partial H_1}{\partial \nu} \right \vert_{\partial M \times (0,T]} =
\left .  \frac{\partial {H}_2}{\partial \nu}  \right \vert_{\partial M \times (0,T]} =\left . \psi \right \vert_{\partial M \times (0,T]}.
\]
Define $F(x,\tau) \doteqdot H_1(x,\tau) - {H}_2(x,\tau)$.
Then $F$ satisfies (\ref{eq F difference unique}) and the boundary condition
$\left .  \frac{\partial F}{\partial \nu} \right
\vert_{\partial M \times (0,T]} = 0$.
Let $\varphi: M \rightarrow \mathbb{R}$ be an arbitrary $C^{2}$ function
which vanishes on $\partial M$.
Fix a $\bar{\tau} \in (0,T]$.
Let $\Phi: M \times\left[ 0,\bar{\tau}\right]  \rightarrow \mathbb{R}$
be the solution to the initial-boundary value problem
(\ref{eq Phi sol backward eq test}) with oblique derivative boundary condition
\[
 \left .  \left ( \frac{\partial \Phi}{\partial \nu} -\langle X, \nu \rangle
 \Phi \right ) \right \vert_{ \partial M \times [0,T]} =0.
\]
The solution always exists.
From (\ref{eq int by parts F}) and similar calculation and argument as in (i)
we conclude that
for $\varepsilon\in\left( 0,\bar{\tau}\right)$
\[
\int_{M }\Phi Fd\mu_{g(\bar{\tau})}= \int_{M}\Phi Fd\mu_{g(\varepsilon)}
\]
and hence $F=0$. This proves ${H}_1 =H_2$.

%%%%%%%%%%%%%%%%%%%%%%%%%%%%%%%%%%%%%%
\vskip .1cm
(iv) Let $\varphi: M \rightarrow\mathbb{R}$ be an arbitrary
 $C^{2}$ function with compact support.
 Fix a $\bar{\tau} \in (0,T]$.
Let $\Phi: M \times\left[ 0,\bar{\tau}\right]  \rightarrow \mathbb{R}$
be the  bounded solution to the initial value problem
(\ref{eq Phi sol backward eq test}).
The solution always exists.

Let $H_i, \, i=1,2$, be two fundamental solution of (\ref{eq fund sol})
on $M \times (0,T]$ centered at $x_0 \in \operatorname{Int}(M)$.
Assume that each $H_i$ satisfies the following equalities:
for any $\Phi$ defined above

\noindent (A1) $\frac{d}{d \tau} \int_M H_i \Phi d \mu_{g(\tau)} = \int_M
\frac{\partial }{\partial \tau} \left (  H_i \Phi d \mu_{g(\tau)} \right )$.

\noindent(A2) $\int_M \Phi \Delta_{g(\tau)} H_i   d \mu_{g(\tau)} =
\int_M H_i  \Delta_{g(\tau)}  \Phi d \mu_{g(\tau)}$.

\noindent (A3) $\int_M  \Phi \nabla_X H_i  d \mu_{g(\tau)} =  \int_M \left (- H_i \nabla_X  \Phi -  H_i\Phi  \operatorname{div}_{g(\tau)} X
 \right )d \mu_{g(\tau)} $.

Define $F(x,\tau) \doteqdot H_1(x,\tau) -{H}_2(x,\tau)$.
Then $F$ satisfies (\ref{eq F difference unique}).
The above assumptions enable us to perform the similar calculation and
argument as in (i), we conclude that
for $\varepsilon\in\left( 0,\bar{\tau}\right)$
\[
\int_{M }\Phi Fd\mu_{g(\bar{\tau})}= \int_{M}\Phi Fd\mu_{g(\varepsilon)}
\]
and hence $F=0$. This proves $H_1 =H_2$.
\end{proof}

\vskip .1cm
\noindent \textbf{Remark} Here we will not discuss when the assumptions
(A1), (A2), and (A3) will be satisfied. The interested readers may
find various estimates of $\Phi$ and $H_i$ from literature which
guarantee that the assumptions hold (see \cite{Ko}, for example).

%%%%%%%%%%%%%%%%%%%%%%%%%%%%%%%%%%%%%%%%%%%%%%%%%%%%%%%%%%%%%%%%%%%%%%%%%%%%
%%%%%%%%%%%%%%%%%%%%%%%%%%%%%%%%%%%%%%%%%%%%%%%%%%%%%%%%%%%%%%%%%%%%%%%%%%%%
%%% \newpage
\section{A local integral estimate of fundamental solutions}

\noindent \textbf{5.1 The cut-off function $h$}.
The following construction of $h$ is adapted from the cutoff function
used by Perelman in his proof of pseudo-locality theorem(\cite[p.25]{Pe02I}),
a similar function is also used by Perelman in his proof of the localized no local
collapsing theorem(\cite[p.21]{Pe02I}).
Let $M^n$, $\Omega$, $g(\tau )$ with $\tau \in \left[ 0,T \right]$,
 $Q$, and $x_0 \in \operatorname{Int}(M)$, be defined as at the beginning of \S2.1.
We also adopt the notations used there.
Let $\tilde{g}$ be a smooth metric on $\Omega$ satisfying
(\ref{eq g tilde g0 equ}) and
$\operatorname{Rc}(\tilde{g}) \geq -K$ on $\Omega$.
 Fix a $r_\ast >0$ such that $\overline{B_{\tilde{g}}(x_0,r_\ast) }
 \subset \Omega$ is compact. We assume
\begin{align}
& \sup_{B_{\tilde{g}}(x_0, r_\ast) \times[0,T]} \vert \mathcal{R}_{ij} (x, \tau)
\vert_{g(\tau)} \leq K_\ast, \label{eq mathcal R bound} \\
& \sup_{B_{\tilde{g}}(x_0, r_\ast) \times[0,T]} \vert \operatorname{Rc}
(x, \tau) \vert_{g(\tau)} \leq K_\ast. \label{eq rc g tau bdd}
\end{align}
By a simple calculation we have
\[
e^{-2K_\ast T} C_0^{-1} \tilde{g} (x) \leq g(x, \tau)
\leq e^{2K_\ast T} C_0 \tilde{g}(x)
\]
for all $x \in B_{\tilde{g}}(x_0, r_\ast)$ and $\tau \in [0,T]$.
Let
\begin{equation}
\hat{r} \doteqdot \frac{1}{2}e^{-2K_\ast T} C_0^{-1} r_\ast. \label{eq r hat def}
\end{equation}
Then closed ball
\[
\overline{B_{g(\tau) }(x_0, \hat{r})} \subset  B_{\tilde{g}}(x_0, r_\ast) \,\,\,\,
\text{ for all } \tau \in [0,T].
\]

For any $x \in B_{g(\tau)}(x_0, \hat{r})$,
we have (see Lemma 18.1 in \cite{CCIII}, for example)
\begin{equation*}
\frac{\partial  }{\partial \tau }
d_{g\left( \tau \right) }\left( x,x_{0}\right)
=\int_{0}^{d_{g\left( \tau \right) }\left( x,x_{0}\right) }
\mathcal{R}_{ij} (\gamma(s), \tau)
\frac{ d \gamma^i}{d s} \frac{ d \gamma^j }{ds} ds.
\end{equation*}
where $\gamma \left( s\right) $ is some unit speed minimal geodesic
with respect to metric $g( \tau ) $ between $x$ and $x_{0}$.
Hence
\begin{equation}
\frac{\partial  }{\partial \tau }
d_{g\left( \tau \right) }\left( x,x_{0}\right) \leq K_\ast d_{g\left( \tau \right) }\left( x,x_{0}\right)
\leq K_\ast \hat{r}. \label{eq deriv dist g tau up est}
\end{equation}

Let $\xi:[0,d_{g(\tau)}(x,x_0)] \rightarrow [0,1]$ be a continuous piecewise smooth function with $\xi(0)=0$ and $\xi(d_{g(\tau)}(x,x_0)) =1$.
We have (see Lemma 18.6 in \cite{CCIII}, for example)
\[
 \Delta _{g(\tau)}d_{g(\tau)}(x,x_0) \leq \int_0^{d_{g(\tau)}(x,x_0) }
 (n-1) \left ( \xi^\prime(s) \right )^2 - \xi^2(s) \operatorname{Rc}
 (\gamma^\prime(s), \gamma^\prime(s) ) ds
\]
Let
\[
x \in  B_{g(\tau)}(x_0, \hat{r}) \setminus  B_{g(\tau)}(x_0,
\frac{1}{10}\hat{r}).
\]
We choose
\begin{equation*}
\xi \left( s\right) =\left\{
\begin{tabular}{ll}
$\frac{10 s}{\hat{r}}$ & $\text{if } s\in [0 ,\frac{\hat{r}}{10}],$ \\
$1$ & if$\text{ }s\in ( \frac{\hat{r}}{10}, d_{g(\tau)}(x,x_0)],$
\end{tabular}
\ \ \right.
\end{equation*}
By a simple calculation and (\ref{eq rc g tau bdd}) we have
\begin{equation}
 \Delta _{g(\tau)}d_{g(\tau)}(x,x_0)  \leq
 \frac{10(n-1)}{\hat{r}} +K_\ast \hat{r} . \label{eq Laplace dist upper bdd}
\end{equation}

Combining (\ref{eq deriv dist g tau up est}) and (\ref{eq Laplace dist upper bdd})
we have proved
\begin{lemma}
Under the assumption given at the beginning of this subsection,
we have
\begin{equation}
\left( \frac{\partial }{\partial \tau }+\Delta _{g\left( \tau \right)
}\right) d_{g\left( \tau \right) }\left( x,x_{0}\right) \leq
 \frac{10(n-1)}{\hat{r}} + 2K_\ast \hat{r}
\end{equation}
for any $x \in  B_{g(\tau)}(x_0, \hat{r})
\setminus  B_{ g(\tau)}(x_0, \frac{1}{10}\hat{r})$.
\end{lemma}

Let $\phi :\mathbb{R}\rightarrow \left[ 0,1\right] $
be a smooth function which is strictly decreasing on
the interval $\left[ 1,2\right] $
and which satisfies
\begin{equation}
\phi \left( s\right) =\left\{
\begin{tabular}{ll}
$1$ & if$\text{ }s\in (-\infty ,1] $ \\
$0$ & if$\text{ }s\in \lbrack 2,\infty )$
\end{tabular}
\ \ \right.   \label{PhiCutoffInfrared01 add1}
\end{equation}
and
\begin{align}
& \left( \phi ^{\prime }(s)\right) ^{2}  \leq  10\phi (s)
\label{eq  temp name 113} \\
& \phi ^{\prime \prime }(s) \geq  -10\phi (s)  \label{eq temp name 114}
\end{align}
for $s\in \mathbb{R}$.

Let $T_1 \in (0,T]$ be a constant to be chosen later (see (\ref{eq T1 choice})).
We define a function $h: \Omega \times [0,T_1] \rightarrow [0,1]$ by

\begin{equation}
h\left( x,\tau \right) =\phi \left( \frac{d_{g\left( \tau \right) }\left(
x,x_{0}\right) +a\left(T_1 -\tau \right) }{b}\right) \label{eq h def}
\end{equation}
where $a$ is a positive constant to be chosen later
(see (\ref{eq a choice})) and $b= \frac{1}{2} \hat{r}$.
Note that $\operatorname{supp}h\left( \cdot , \tau \right) \subset
B_{g\left( \tau \right) }\left( x_{0},\hat{r} \right) $.
Let $X(\tau)$, $\tau \in [0,T]$, be a smooth family of vector fields on $B_{\tilde{g}}(x_0, r_\ast)$, and let
\[
w\left( x,\tau \right) \doteqdot \frac{d_{g\left(
\tau \right) }\left( x,x_{0}\right) +a\left(T_1 -\tau \right) }{b}.
\]
We compute
\begin{eqnarray}
\frac{\partial h}{\partial \tau }+\Delta _{g\left( \tau \right) }h -\nabla_X h
&=&\frac{\phi ^{\prime }\left( w\right) }{b}\left( \left(
\frac{\partial }{\partial\tau }
+\Delta _{g\left( \tau \right) } -\nabla_X \right) d_{g\left( \tau \right)
}\left( x,x_{0}\right) -a\right)  \notag  \\
&&+\frac{\phi^{\prime \prime}}{b^{2}}\left\vert \nabla _{
g\left( \tau \right) }d_{g\left(
\tau \right) }\left( x,x_{0}\right) \right\vert _{g\left( \tau \right) }^{2}.
\label{eq tem use 115}
\end{eqnarray}
First we choose
\begin{equation}
a \geq  \frac{10(n-1)}{\hat{r}} + 2K_\ast \hat{r} +K_1 \label{eq a choice}
\end{equation}
where $K_1 \doteqdot \sup_{B_{\tilde{g}}(x_0, r_\ast) \times [0,T]}
|X|_{g(\tau)}$.
Next we choose $T_1$ such that
\begin{equation}
\frac{1}{10}\hat{r} +aT_1 \leq b= \frac{1}{2} \hat{r}. \label{eq T1 choice}
\end{equation}
Then $w(x,\tau) \leq 1$ and $\phi^\prime (w) =0$ for any
$x \in B_{ g(\tau)}(x_0, \frac{1}{10}\hat{r})$ and $\tau \in [0,T_1]$.
We have proved that either  $\phi^\prime (w(x, \tau)) =0$ or
\[
\left( \frac{\partial }{\partial
\tau }+\Delta _{g\left( \tau \right) } -\nabla_X \right) d_{g\left( \tau \right)
}\left( x,x_{0}\right) - a \leq 0,
\]
i.e.,
 \begin{equation}
 \frac{\phi ^{\prime }\left( w\right) }{b}\left( \left(
\frac{\partial }{\partial
\tau }+\Delta _{g\left( \tau \right) }- \nabla_X \right) d_{g\left( \tau \right)
}\left( x,x_{0}\right) -a\right) \geq 0 \label{eq dist heat oper act}
 \end{equation}
 for $x \in M$ and $\tau \in [0,T_1]$.

Combining (\ref{eq tem use 115}), (\ref{eq dist heat oper act})  and
$\left\vert \nabla _{g\left( \tau \right) }d_{g\left( \tau \right)
}\left( x,x_{0}\right) \right\vert _{g\left( \tau \right) }=1$
we have proved
\begin{lemma} \label{lem cutoff funct h}
Under the assumption given at the beginning of this subsection,
the function $h$ defined in (\ref{eq h def}) with the choice of $a$,
$b$ and $T_1$ given by (\ref{eq a choice} ) and
(\ref{eq T1 choice}) satisfies
\begin{equation*}
\frac{\partial h}{\partial \tau }+\Delta _{g\left( \tau \right) }h
- \nabla_X h \geq - \frac{10}{b^{2}}h.
\end{equation*}
\end{lemma}

%%%%%%%%%%%%%%%%%%%%%%%%%%%%%%%%%%%%%%%%%%%%%%%
%%%%\vskip .3cm
\noindent \textbf{5.2 Lower bound of integrals of fundament solutions on balls}.

In this subsection we use the setup described at the beginning of
\S 3.
 Let $r_\ast$ in \S 5.1 be a positive constant such that $\overline{B_{\tilde{g}}
(x_0, r_\ast)}$ is compact subset of $M$ and such that
 (\ref{eq mathcal R bound})  and (\ref{eq rc g tau bdd}) hold.
Let $h(x,\tau) = \phi(w(x,\tau))$ be the function defined by (\ref{eq h def})
with structure constants $\hat{r}$ defined by (\ref{eq r hat def}),
$a$ satisfying (\ref{eq a choice}), $b=\frac{1}{2} \hat{r}$,
and $T_1$ satisfying (\ref{eq T1 choice}).
Note that the support $\operatorname{supp} h(\cdot, \tau) \subset \overline{B_{\tilde{g}}(x_{0}, r_\ast)}$ for each $\tau \in [0,T_1]$.

Define $K_2$ by
\begin{equation}
K_2 \doteqdot \sup_{  B_{\tilde{g}}(x_{0},
r_\ast )  \times \left[ 0,T_1 \right] } \left \{ Q,
\operatorname{div}_{g(\tau)} X \right \}. \label{eq K 2 def}
\end{equation}
Let $u$ be a nonnegative solution of (\ref{eq heat gen}) on $M \times (0,T]$.
We compute
\begin{align*}
\frac{d}{d\tau}\int_{M} u h d\mu_{g\left( \tau\right)  }
=&\int_{M}h \left(   \frac{\partial u}{\partial\tau}-\Delta_{g\left(
\tau\right)  } u+ \nabla_X u +Qu \right ) d\mu_{g\left( \tau\right)  }  \\
& +\int_{M} u \left(
\begin{array}
[c]{c}
\frac{\partial h}{\partial \tau}+\Delta_{g\left(\tau\right) }h -\nabla_X h \\
+\left(  \mathcal{R} -Q -\operatorname{div}_{g(\tau)} X \right)  h
\end{array}
\right)d\mu_{g \left(  \tau\right)  }  \\
\geq & -\left ( \frac{10}{b^2} + n K_\ast + K_2 \right ) \int_{M}
u hd\mu_{g\left( \tau\right)  },
\end{align*}
where  we have used (\ref{eq mathcal R bound}) and Lemma  \ref{lem cutoff funct h}
to get the last inequality.
We have proved
\[
\frac{d}{d\tau }\ln \int_{M} u h d\mu_{g\left( \tau \right) }
\geq   -\left ( \frac{10}{b^2} + n K_\ast + K_2 \right ).
\]
Hence for any $0 < \tau_1 < \tau \leq T_1$ we have
\begin{equation}
 \int_{M } u h d\mu _{ g \left( \tau \right) }
 \geq e^{- \left ( \frac{10}{b^2} + n K_\ast + K_2 \right ) \left( \tau -\tau _{1}\right) } \int_{M } u h d\mu_{
g \left( \tau _{1}\right) }. \label{eq uk hk int ineq}
\end{equation}

Suppose that $u$ can be continuously extended to a function defined
on $M \times [0,T]$,
then we have
\begin{equation*}
\int_{M } u h d\mu _{ g \left( \tau \right) }
\geq e^{- \left ( \frac{10}{b^2} + n K_\ast + K_2 \right ) T_1 }
\int_{ M } u h  d\mu_{ g \left( 0\right) }
\text{~~~~ for } \tau \in (0,T_1].
\end{equation*}
Since $\operatorname{supp}h (\cdot,\tau) \subset B_{\tilde{g}}(x_{0},
 r_\ast)$ for any $\tau \in [0,T_1]$,
 using (\ref{eq T1 choice}) we have proved the following.
 \begin{proposition} \label{prop local int fund sol lower bdd}
Let $g(\tau), \, \tau \in [0,T ]$, be a smooth family of smooth metrics
on $M^n$ and let $\tilde{g}$ be a smooth metric on $M$.
We assume that $\overline{B_{\tilde{g}} (x_0, r_\ast)}$ is compact subset of $M$.
Let $u$ be a nonnegative solution of (\ref{eq heat gen}) on $M \times [0,T]$.
Then
\begin{equation*}
\int_{B_{\tilde{g} }\left( x_{0}, r_\ast \right) }
 u d\mu _{ g \left( \tau \right) } \geq e^{- \left (
 \frac{10}{b^2} + n K_\ast + K_2 \right ) T_1 }
 \int_{B_{\tilde{g}}(x_{0}, \frac{1}{10} \hat{r} )}u d\mu_{ g \left( 0\right) }
\text{~~~~ for } \tau \in (0,T_1].
\end{equation*}
Here $b=\frac{1}{2} \hat{r}$ is defined by (\ref{eq r hat def}),
$K_*$ is defined by (\ref{eq mathcal R bound}) and (\ref{eq rc g tau bdd}),
$K_2$ is defined by (\ref {eq K 2 def}),
and $T_1$ is defined by (\ref{eq a choice}) and (\ref{eq T1 choice}).
\end{proposition}

When $H$ is a fundamental solution of
(\ref{eq fund sol}) on $M \times (0,T]$ centered at $x_0 \in \operatorname{Int}(M)$,
we have
\[
\lim_{\tau _{1} \rightarrow 0} \int_{M } H (x, \tau_1) h (x,\tau_1)
 d\mu_{ g \left( \tau _{1}\right) } (x) =h (x_0, 0) =1.
\]
Then it follows from (\ref{eq uk hk int ineq}) that
\begin{equation*}
\int_{M } H h d\mu _{ g \left( \tau \right) }
\geq e^{- \left ( \frac{10}{b^2} + n K_\ast + K_2 \right ) T_1 }
\text{~~~~ for } \tau \in (0,T_1].
\end{equation*}
Hence we have proved the following.

 \begin{corollary} \label{cor local int fund sol lower bdd}
Let $g(\tau), \, \tau \in [0,T]$, be a smooth family of smooth metrics
on $M^n$ and let $\tilde{g}$ be a smooth metric on $M$.
We assume that $\overline{B_{\tilde{g}} (x_0, r_\ast)}$ is compact subset in $M$.
Let $H$ be a fundamental solution of (\ref{eq fund sol}) on $M \times (0,T]$
centered at $x_0 \in \operatorname{Int}(M)$. Then
\begin{equation*}
\int_{B_{\tilde{g} }\left( x_{0}, r_\ast \right) }
 H d\mu _{ g \left( \tau \right) } \geq e^{- \left (
 \frac{10}{b^2} + n K_\ast + K_2 \right ) T_1 }
\text{~~~~ for } \tau \in ( 0,T_1].
\end{equation*}
Here $b=\frac{1}{2} \hat{r}$ is defined by (\ref{eq r hat def}),
$K_*$ is defined by (\ref{eq mathcal R bound}) and (\ref{eq rc g tau bdd}),
$K_2$ is defined by (\ref {eq K 2 def}),
and $T_1$ is defined by (\ref{eq a choice}) and (\ref{eq T1 choice}).
\end{corollary}

%%%%%%%%%%%%%%%%%%%%%%%%%%%%%%%%%%%%%%%%%%%%%%%%%%%%%%%%%%%%%%%%%%%%%%%%%%%%%%%
\section*{Acknowledgments}
{Part of this work was done while the author was visiting Department of Mathematics,
University of California at San Diego in early 2009. The author thanks
Bennett Chow and Lei Ni for their invitation and hospitality and
thanks Bennett Chow for the helpful discussion.}

%%%%%%%%%%%%%%%%%%%%%%%%%%%%%%%%%%%%%%%%%%%%%%%%%%%%%%%%%%%%%%%%%%%%%%%%%%%%%%%
\bibliographystyle{natbib}

%%%%%%%%%%%%%%%%%%%%%%%%%%%%%%%%%%%%%%%%%%%%%%%%%%%%%%%%%%%%%%%%%%%%%%%%%%%%%%%
%%%%%%%%%%%%%%%%%%%%%%%%%%%%%%%%%%%%%%%%%%%%%%%%%%%%%%%%%%%%%%%%%%%%%%%%%%%%%%%

%%%% \begin{center} {\bf Keywords:} {Cheeger-Gromov convergence,
%%%% heat kernels, Ricci flow \hfill}
%%%% \end{center}
%%%%\begin{abstract} In this note we show the convergence of the 
%%%% fundamental solutions of the heat-type equations assuming 
%%%% the Cheeger-Gromove convergence of the underlying manifolds  
%%%% and the uniform $L^1$-bound. We also prove a local
%%% integral estimate of fundamental solutions. 
%%%%\end{abstract}

\medskip
\noindent {\sc Peng Lu}, University of Oregon

\noindent e-mail: penglu@uoregon.edu

\end{document}